\documentclass[preview,11pt]{article}

\usepackage[nottoc]{tocbibind}

\usepackage{amssymb,latexsym}
\usepackage{graphicx}
\usepackage{amsfonts, amsmath, amssymb, dsfont, setspace}
\usepackage{fixltx2e}

\usepackage{color}
\usepackage[dvipsnames]{xcolor}
\usepackage{verbatim}
\usepackage{vmargin}
\usepackage{url}
\usepackage[english]{babel}
\usepackage[latin1]{inputenc}
\usepackage[T1]{fontenc}
\usepackage{enumitem}
\usepackage{array}
\usepackage{fancyhdr}
\usepackage{lmodern}
\usepackage{cleveref}
\crefformat{equation}{(#2#1#3)}
\usepackage[font={footnotesize}]{caption}
 \usepackage[normalem]{ulem}
\newtheorem{Assumption}{Assumption}[section]
\newtheorem{remark}{Remark}[section]

\usepackage{natbib}

\title{Numerical continuation of spiral waves in \\ heteroclinic networks of cyclic dominance}
\author{Cris R. Hasan, \\[-1mm]
  \small{School of Mathematical Sciences, University College Cork, Ireland} \\[-1mm]
  \small{rhasan@ucc.ie}
  \and
  Hinke M. Osinga, Claire M. Postlethwaite \\[-1mm]
  \small{Department of Mathematics, The University of Auckland, New Zealand}
  \and
  and Alastair M. Rucklidge \\[-1mm]
  \small{School of Mathematics, University of Leeds, UK}}
\date{}

\begin{document}

\thispagestyle{empty}

\maketitle

 \begin{abstract}
Heteroclinic-induced spiral waves may arise in systems of partial differential equations that exhibit robust heteroclinic cycles between spatially uniform equilibria. Robust heteroclinic cycles arise naturally in systems with invariant subspaces and their robustness is considered with respect to perturbations that preserve these invariances. We make use of particular symmetries in the system to formulate a relatively low-dimensional spatial two-point boundary-value problem in Fourier space that can be solved efficiently in conjunction with numerical continuation. Our numerical set-up is formulated initially on an annulus with small inner radius, and Neumann boundary conditions are used on both inner and outer radial boundaries. We derive and implement alternative boundary conditions that allow for continuing the inner radius to zero and so compute spiral waves on a full disk. As our primary example, we investigate the formation of heteroclinic-induced spiral waves in a reaction-diffusion model that describes the spatiotemporal evolution of three competing populations in a two-dimensional spatial domain---much like the Rock--Paper--Scissors game. We further illustrate the efficiency of our method with the computation of spiral waves in a larger network of cyclic dominance between five competing species, which describes the so-called Rock--Paper--Scissors--Lizard--Spock game.
 \end{abstract}

\section{Introduction}

In dynamical systems, a heteroclinic cycle is a set of trajectories that
connect equilibria in a topological circle~\citep{guckenheimer1988, krupa1997}.
In general, such a cycle does not persist under perturbation, unless the
dynamical system has a special structure. More precisely, this phenomenon
occurs in systems with special properties that allow for the existence
of a sequence of invariant subspaces. A heteroclinic cycle is called robust, or
structurally stable, when it persists under small perturbations that preserve
this special structure. In one of its simplest forms, a robust heteroclinic
cycle involves three saddle equilibria and their connecting trajectories in a
system for which the equilibria are pairwise contained in an invariant
subspace. This situation often occurs in population models and cyclic interactions are known to provide a naturally selective coexistence
mechanism for interspecific competition between three subpopulations, including
morphs of the side-blotched lizard~\citep{Sinervo1996, Sinervo2000}, coral reef
invertebrates~\citep{Jackson1975, Taylor1990}, and strains of \emph{Escerichia
coli}~\citep{Kerr2002, Kirkup2004}. The dynamics near a robust three-equilibrium heteroclinic cycle also emulates the famous children's game of
Rock--Paper--Scissors, where Rock crushes Scissors, Scissors cut Paper and
Paper wraps Rock.

The focus of this paper is on the computation of spiral waves arising from cyclic interactions between
competing populations in the presence of spatial diffusion. We assume that this
system is modelled by a system of partial differential equations (PDEs) given
in vector form as
 \begin{equation}
  \label{eq:VectorForm}
  \mathbf{U}_t = \mathbf{f}(\mathbf{U}) + \mathbf{U}_{xx}+\mathbf{U}_{yy}.
 \end{equation}
Here, $\mathbf{U}(x,y,t) \in \mathbb{R}^m$ represents $m$ different species,
$\mathbf{f} \colon \mathbb{R}^m \to \mathbb{R}^m$ is a sufficiently smooth function expressing the nonlinear kinetics (i.e., local species
interactions), and the spatial derivatives $\mathbf{U}_{xx}+~\mathbf{U}_{yy}$
are the Laplacian terms modelling diffusion. We are interested in a specific class of \cref{eq:VectorForm} that satisfies the following assumption:
\begin{Assumption} 
\label{assumption:symmetry}
\mbox{} \\
We assume that, in the absence of spatial variation and diffusion terms, \cref{eq:VectorForm} admits an attracting, robust heteroclinic cycle between equilibria solutions. Furthermore, we assume that $\mathbf{f}$  is invariant under cyclic permutations of its arguments.
\end{Assumption}

 Systems of the form \cref{eq:VectorForm} that satisfy \Cref{assumption:symmetry} describe a large class of problems that feature heteroclinic cycles with permutation symmetries (i.e., cyclic symmetries). Well-known examples of this class include those arising in equivariant bifurcation theory (e.g., the three-variable Guckenheimer--Holmes cycle~\citep{guckenheimer1988} and the four-variable Field--Swift cycle~\citep{field1991}); population models (e.g., the May--Leonard model~\citep{May1975} and other higher-dimensional generalisations of Lotka--Volterra type); and also those that can appear in other mathematical modelling problems (e.g., the example of \citet{proctor1988} from fluid mechanics).
Note that the diffusion coefficients for all state variables in system~\cref{eq:VectorForm} are equal; this is necessary to preserve the permutation symmetry. If we perturb the system by breaking this symmetry, for example by choosing unequal diffusion coefficients slightly different from 1, we expect the dynamics to be similar, although the computations will be more expensive.

\begin{remark}
The computational method presented here can be generalized to compute spiral waves
in any reaction--diffusion system of the general form \cref{eq:VectorForm}, but in the absence of the cyclic symmetry
given by \Cref{assumption:symmetry} the computational cost will then likely be larger.
\end{remark}
%
\begin{figure}[t!] 
  \centering
  \includegraphics[scale=1.0]{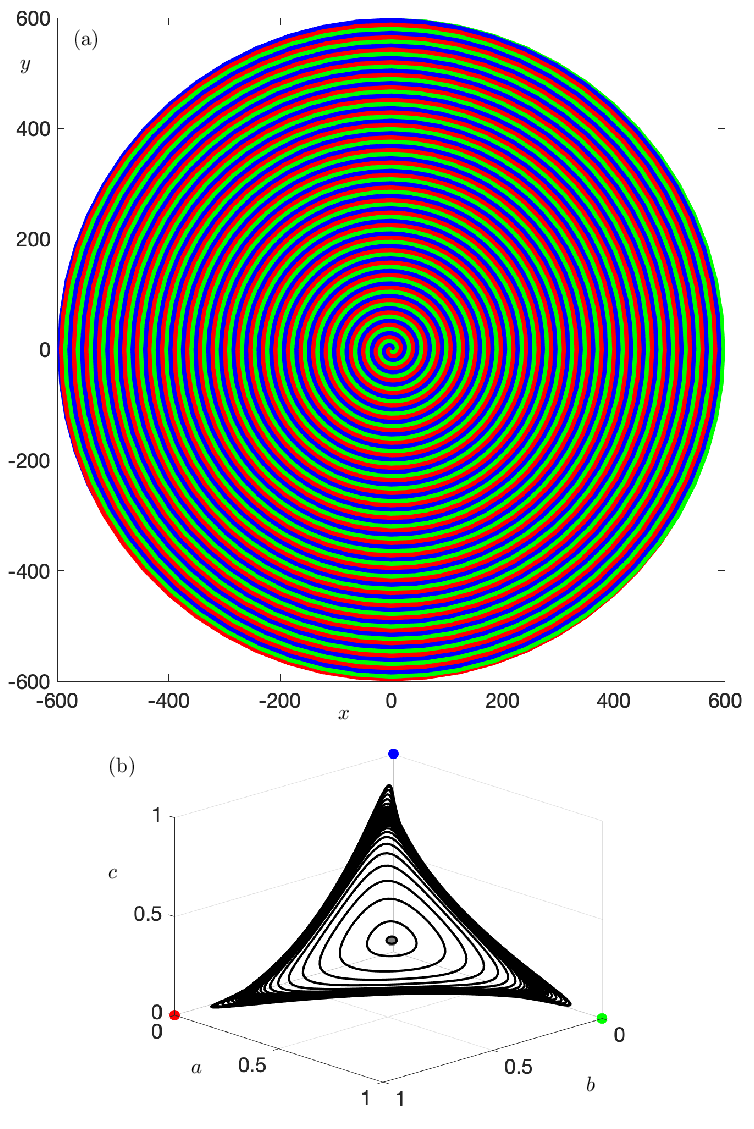}
  \caption{\label{fig:LargeSW}
    A large spiral wave formed by the three competing species of
    system~\cref{eq:PDElaboratory} with parameters $\sigma=3.2$ and
    $\zeta=0.8$, shown at a fixed moment in time. Panel~(a) shows the spiral
    wave solution on an annulus in the ($x,y$)-plane; the three different
    colours represent which of the population densities $a$, $b$ and $c$ is the
    largest. Panel~(b) shows their distribution in $(a, b, c)$-space along a
    selection of concentric circles around the origin in the ($x,y$)-plane; the
    smallest and largest circles are the inner and outer boundaries of the
    computed domain at radii $0.01$ and $600$, respectively.
    The grey dot is the coexistence equilibrium and the colored dots are the on-axis equilibria where only one species survives.}
 \end{figure}

\Cref{fig:LargeSW} shows an example of a spiral wave for the simplest case of system~\cref{eq:VectorForm} with $m = 3$; the precise equations are given in equation \cref{eq:PDElaboratory} in the next section.
The species are represented in terms of their scaled
densities $a$, $b$ and $c$ at an arbitrarily chosen instant in time. With $m =
3$ there exists only one heteroclinic cycle, formed by connecting orbits
between the equilibria $(a, b, c) = (1, 0, 0)$, $(0, 1, 0)$ and $(0, 0, 1)$.
There also exists a spatially-uniform coexistence equilibrium $(a, b, c) =
\frac{1}{3 + \sigma} (1, 1, 1)$, at which all species survive.
The spiral wave was computed on an annulus in the $(x, y)$-plane, centred at the origin, with inner radius $0.01$ and outer radius $600$. Panel~(a)
illustrates the spiral wave in the ($x,y$)-plane with red, green and blue
regions indicating the dominance of $a$, $b$ and $c$, respectively, that is, which of the three populations, $a$, $b$ or $c$ is the largest. 
Panel~(b)
shows the population densities (black curves) in $(a, b, c)$-space, where $x$ and $y$ vary along a family of concentric circles around the origin; here, the smallest cycle in $(a, b, c)$-space corresponds to the inner boundary with radius $0.01$ in the $(x, y)$-plane and the largest corresponds to the outer
boundary of the annulus with the very large radius of $600$. 
The grey dot is the coexistence equilibrium and the coloured dots are the on-axis equilibria that are involved in the limiting heteroclinic cycle.
The periodic curves in the $(a, b, c)$-space are centred on a point which is close to, but not exactly at, the coexistence equilibrium, indicating that the core of the spiral at  $(x,y)=(0,0)$ should satisfy $a = b = c \approx \frac{1}{3 + \sigma}$.
The population densities along the outermost circle accumulate on a periodic solution that comes close to, but remains a finite distance from, each of the three 
equilibria in turn. This periodic solution is related to the periodic 
travelling waves found in the one-dimensional version of this problem, and to 
the underlying heteroclinic cycle \citep{hasan2021spatiotemporal}.
In contrast to the almost uniform
distribution along the inner circle, the variation of the population densities
along the outer circle is non-uniform, with each of the species being
close to zero along two thirds of the perimeter.

Most previous approaches for numerical continuation of
spirals~\citep*{ barkley1992, bar2003,Bordyugov2007,Dodson} have focussed on
systems of reaction--diffusion type, in which a spatially uniform equilibrium
solution loses stability in a Hopf (Turing) bifurcation, leading to spiral waves. The
spiral waves in our examples have oscillations associated with heteroclinic
cycles rather than Hopf bifurcations, so we call these
\emph{heteroclinic-induced spirals}. 
As the radius of the spiral goes to infinity, the system approaches a periodic orbit in the far field that is associated with a heteroclinic bifurcation from the heteroclinic cycle~\citep{ postlethwaite2017, postlethwaite2019,hasan2021spatiotemporal}.
The proximity to a heteroclinic cycle makes the continuation the spiral wave a challenge
because in the far field, the variables switch fairly abruptly from being
nearly zero to being nearly one and back again.

Our method is based on the approach by \cite{Bordyugov2007}, who compute spiral waves in reaction--diffusion
systems as a continuation of solutions to a two-point boundary-value problem
(BVP) formulated in Fourier space. 
In systems of heteroclinic networks, when the periodic orbit in the far field approaches the heteroclinic cycle,
the rapid switching between episodes of the variables being nearly
constant means that a large number of Fourier modes may be necessary. Consequently,
the computation of the spiral wave on a large domain requires solving a
computationally expensive~\hbox{BVP}.
We make use of the symmetry of the heteroclinic
cycle, imposing phase relations between the Fourier modes, and so obtaining a
significant reduction in the dimension of the~\hbox{BVP}; the efficiency of this approach is most apparent when a large number of species is considered.

We also aim to study the dynamics on the full disk rather than just on an annulus, which is not possible using the Neumann boundary conditions suggested by \citet{Bordyugov2007}.
More precisely, when taking the inner boundary to zero radius, one encounters a singular Laplacian with such boundary conditions.
We derive and implement the correct boundary conditions at the core in order to obtain a bounded Laplacian term in polar coordinates.
In this way, we compute the spiral waves on a full disk, as opposed to considering an annulus with a small hole around the origin. 

As our primary example, we use the three-species model illustrated in
\cref{fig:LargeSW} and introduced in the next section. This is the simplest
model with a robust heteroclinic cycle and the same model that we studied
in \citet{postlethwaite2017, postlethwaite2019} and \citet{hasan2021spatiotemporal}, where we focussed on
the bifurcation structure of different types of heteroclinic cycles, and the
computation of periodic travelling waves. The dynamics of the
heteroclinic-induced spiral waves for this model have also been investigated
by others \citep{reichenbach2008, Szczesny2013, Szczesny2014, Szolnoki2014}, with
focus on existence and break-up of the spiral-wave patterns in the presence of small defects, based on analysis and simulation of the full~\hbox{PDE}. We
chose this three-species competition model as our primary example because the
numerical set-up and efficiency gains are already evident for this simple
example, and the computations are readily compared with other results from the literature. We also show how to adapt the numerical set-up for the computation of spiral waves in a system of the form~\cref{eq:VectorForm} with $m = 5$; see already equation \cref{eq:5speciesPDE}. This higher-dimensional example not only has a heteroclinic cycle between all five equilibria with only one surviving species, but also one between five equilibria with three surviving species each. We compute
spiral waves from both families to emphasise the versatility of our computational
method.

This paper is organised as follows. In the next section, we introduce our
primary example given by system~\cref{eq:VectorForm} with $m = 3$ and review the
method from \citet{Bordyugov2007}. 
In \cref{sec:symmetry}, we exploit the
cyclic symmetry and present a reduced~\hbox{BVP}. 
We discuss the boundary conditions at the core of the spiral waves in \cref{sec:BCs}.
In \cref{sec:comparethreespecies}, we present a case study
where we explore the properties of spiral waves for our primary example and compare with published results.
In \cref{sec:5species}, we introduce the system of associated competing population model with cyclic dominance between five
species and implement the modified continuation method to compute two families
of five-component spiral waves. We conclude the paper with a discussion and
final remarks in \cref{sec:discussion}.
\section{Computing spiral waves in the three-species model}
\label{sec:background}
The model for heteroclinic-induced spiral waves between three competing
populations was first proposed by  \cite{reichenbach2007} and is based
on the system of three ordinary differential equations introduced by \citet{May1975} as a
model of competing populations without spatial structure and 
diffusion. See~\citep{Frey2010, Szolnoki2014} for details and a
recent review. The system of partial differential equations (PDEs) is defined
in terms of three population densities $a(x,y,t), b(x,y,t), c(x,y,t) \geq 0$
that are scaled to unity; it is given by
 \begin{equation}
  \label{eq:PDElaboratory} 
  \left\{
    \begin{array}{rrlll}
      \dot{a} &=&  a \, (1 - a - b - c - (\sigma+\zeta) \, b + \zeta \, c) + \nabla^2 a, \\
      \dot{b} &=&  b \, (1 - a - b - c - (\sigma+\zeta) \, c + \zeta \, a) + \nabla^2 b, \\
      \dot{c} &=&  c \, (1 - a - b - c - (\sigma+\zeta) \, a + \zeta \, b) + \nabla^2 c.
    \end{array} \right.
 \end{equation}
 \noindent
 The parameters $\sigma$ and $\zeta$ are non-negative constants that represent
 removal and replacement rates of two different interacting species,
 respectively \citep{Szczesny2013}.
 We used the notation $\nabla^2$ for the Laplace operator
 ($\frac{\partial^2}{\partial x^2} + \frac{\partial^2}{\partial y^2}$) that
 models the diffusion on the two-dimensional domain; nonlinearity in the
 diffusion \citep{Szczesny2013}, is not included in this model. The only nonlinearity is the
 population kinetics, as defined by the function $\mathbf{f}$ in
 system~\cref{eq:VectorForm}. Note the cyclic symmetry of the kinetic term:
 if $(a, b, c) = (d_1, d_2, d_3)$ is a solution to
 system~\cref{eq:PDElaboratory} then the cyclic permutations $(a, b, c) = (d_2,
 d_3, d_1)$ and $(a, b, c) = (d_3, d_1, d_2)$ are also solutions to this system
 of PDEs. Furthermore, in the absence of spatial distribution, the coexistence 
equilibrium point is unstable and all trajectories are attracted to the 
heteroclinic cycle. Subspaces with one or more population equal 
to zero are invariant, and hence,
 this system is of the form~\cref{eq:VectorForm} and satisfies  \Cref{assumption:symmetry}.

Spiral waves for system~\cref{eq:PDElaboratory} can be
computed with any of the methods described by \cite{barkley1992}, ~\citet*{bar2003}, \cite{Bordyugov2007} or \cite{Dodson}. In this section, we only review the numerical continuation method introduced by \cite{Bordyugov2007}, because it forms the basis for our improved approach.

\subsection{Spiral waves on an annulus: the Fourier decomposition}
\label{sec:overview}

The spiral waves for systems of the form~\cref{eq:VectorForm} are so-called
rigidly rotating spirals, which are periodic in time and any shift in time is
equivalent to a rotation in space. Furthermore, the centre of rotation, that
is, the tip of the spiral, can be anywhere in the $(x, y)$-plane. To
capture the spatiotemporal rotation symmetry, it is convenient to re-write
system~\cref{eq:VectorForm} in polar coordinates, where we define $x = r
\cos{\phi}$ and $y = r \sin{\phi}$, with $r=0$ being the tip of the spiral. Then the PDE is given by
 \begin{displaymath}
  \mathbf{U}_t = \mathbf{f}(\mathbf{U}) + \mathbf{U}_{rr} +\frac{1}{r} \mathbf{U}_{r} + \frac{1}{r^2} \mathbf{U}_{\phi \phi}.
 \end{displaymath}
The next natural step is to assume that the spiral waves rotate at a constant angular frequency $\omega$ and introduce the co-rotating variable $\theta = \phi + \omega \, t$. Then $\frac{\partial}{\partial \phi} \mapsto \frac{\partial}{\partial \theta}$ and $\frac{\partial}{\partial t} \mapsto \omega \frac{\partial}{\partial \theta} + \frac{\partial}{\partial t}$, which leads to the PDE in the co-rotating frame of reference:
 \begin{displaymath}
  \mathbf{U}_t = \mathbf{f}(\mathbf{U}) + \mathbf{U}_{rr} + \frac{1}{r} \mathbf{U}_{r} +
                           \frac{1}{r^2} \mathbf{U}_{\theta \theta} - \omega \mathbf{U}_{\theta}.
 \end{displaymath} 
A rigid spiral wave with angular frequency $\omega$ and temporal period $T=2\pi/\omega$ is a stationary solution of
this PDE, i.e., it satisfies
 \begin{equation}
  \label{eq:stationary}
  \mathbf{f}(\mathbf{U}) + \mathbf{U}_{rr} + \frac{1}{r} \mathbf{U}_{r} +
  \frac{1}{r^2} \mathbf{U}_{\theta \theta} - \omega \mathbf{U}_{\theta} = 0.
 \end{equation} 
Since spiral wave solutions are periodic in $\theta$, that is, $\mathbf{U}(r,
\theta) = \mathbf{U}(r, \theta + 2 \pi)$, they can be expressed as a Fourier
series expansion with respect to the second argument. We assume that the
Fourier coefficients converge rapidly so that only a finite number of modes is
sufficient to approximate $\mathbf{U}(r, \theta)$. Using $N$ Fourier modes,
where $N$ is assumed to be even, we define $N$ uniformly spaced angles $\theta
= 2\pi \frac{n}{N}$, with $n = 0, 1, \dots, N-1$, and approximate the spiral
wave as
 \begin{displaymath} 
  \mathbf{U}(r, 2\pi \tfrac{n}{N}) = \sum\limits_{k =0}^{N-1} \mathbf{\widehat{U}}(r, k)  \, e^{2\pi \, i \frac{k \, n}{N}}.
 \end{displaymath} 
Here, each $\mathbf{\widehat{U}}(r, k)$, with $k = 0, 1, \dots, N-1$, is a vector of (complex-valued) Fourier coefficients associated with the $k$th Fourier mode, and it is defined as
\begin{displaymath} 
  \mathbf{\widehat{U}}(r, k) = \frac{1}{N} \sum\limits_{n=0}^{N-1} \mathbf{U}(r, 2\pi \tfrac{n}{N}) \, e^{-2\pi \, i \frac{k \, n}{N}}.
\end{displaymath} 
The goal is to solve the PDE in the co-rotating frame~\cref{eq:stationary} with respect to the unknown Fourier coefficient vectors $\mathbf{\widehat{U}}(r, k)$ and the unkown angular frequency $\omega$.  To this end, we need the Fourier transform of the nonlinear kinetics function $\mathbf{f}$. For the case of the three-species population model~\cref{eq:PDElaboratory}, the nonlinear kinetics is quadratic and Fourier coefficients $\mathbf{\widehat{f}}(r,k)$, for $k = 0, 1, \dots, N-1$, of the discretised Fourier approximation
\begin{displaymath}
  \mathbf{f\Big(U}(r, 2\pi \tfrac{n}{N})\Big) = \sum\limits_{k=0}^{N-1} \mathbf{\widehat{f}}(r, k) \, e^{2\pi \, i \frac{k \, n}{N}}
\end{displaymath}
can formally be expressed as a sum of convolutions in terms of the Fourier transforms $\widehat{a}$, $\widehat{b}$ and $\widehat{c}$ of the population densities $a$, $b$ and $c$, respectively. More precisely,
\begin{equation}
  \label{eq:nonlinearTerm}
  \mathbf{\widehat{f}}(r, k) = \begin{pmatrix}
    \widehat{a} - \widehat{a} \ast \widehat{a} - (1 + \sigma +\zeta) \, \widehat{a} \ast \widehat{b} - (1-\zeta) \,  \widehat{a} \ast \widehat{c} \\
    \widehat{b} - \widehat{b} \ast \widehat{b} - (1 + \sigma +\zeta) \, \widehat{b} \ast \widehat{c} - (1-\zeta) \, \widehat{b} \ast \widehat{a} \\
    \widehat{c} - \widehat{c} \ast \widehat{c} - (1 + \sigma+\zeta)\, \widehat{c} \ast \widehat{a} - (1-\zeta) \, \widehat{c} \ast \widehat{b}
  \end{pmatrix} (r,k),
\end{equation}
where, for example, the convolution $\widehat{a}(r,k) \ast \widehat{b}(r,k)$ is calculated as a fast Fourier transform (FFT) of the product $a(r,\theta) \, b(r,\theta)$ obtained from the inverse fast Fourier transform (IFFT) of the individual Fourier coefficients $\widehat{a}(r,k)$ and $\widehat{b}(r,k)$; in other words,
\begin{displaymath}
  \widehat{a}(r,k) \ast \widehat{b}(r,k) = {\rm FFT} \left( {\rm IFFT} \left[ \widehat{a}(r,k) \right] \,
                                                               {\rm IFFT} \left[ \widehat{b}(r,k) \right] \right).
\end{displaymath}
Hence, in Fourier space, the PDE~\cref{eq:stationary} is given by the large system of $m \, N$ second-order complex equations
\begin{displaymath}
  \mathbf{\widehat{f}}(r,k) + \mathbf{\widehat{U}}_{rr}(r,k) + \frac{1}{r} \, \mathbf{\widehat{U}}_{r}(r,k)
                - \frac{k^2}{r^2} \, \mathbf{\widehat{U}}(r,k) - i \,  k \, \omega \, \mathbf{\widehat{U}}(r,k) = 0,
\end{displaymath}
for the $N$ unknown Fourier coefficient vectors $\mathbf{\widehat{U}}(r,k) \in \mathbb{C}^m$. Note that only derivatives with respect to the radial coordinate remain. Hence, we can express the PDE~\cref{eq:stationary} as a system of first-order ordinary differential equations (ODEs):
\begin{equation}
  \label{eq:autonomousODE}
  \left\{ \begin{array}{llll}
	    \mathbf{\widehat{U}}^{\prime}    &=&  \mathbf{\widehat{U}}_r,	 \\
            \mathbf{\widehat{U}}_r^{\prime} &=&  -\mathbf{\widehat{f}} - \frac{1}{r} \, \mathbf{\widehat{U}}_r
                                                  + \frac{k^2}{r^2} \, \mathbf{\widehat{U}} + i \, k \, \omega \, \mathbf{\widehat{U}}, \\
            r^{\prime}                                    &=&  1,
          \end{array} \right.
\end{equation}
where $^\prime$ represents derivation with respect to $r$ and the last equation renders the system autonomous. When split into real and imaginary parts, we now have $4 \, m \, N + 1$ equations---$2 \, m$ for each of the real and imaginary parts of the $N$ Fourier modes---and $4 \, m \, N + 1$ unknowns---the real and imaginary parts of the $m$-dimensional Fourier coefficient vectors $\mathbf{\widehat{U}}(r,k)$, their $m$-dimensional derivatives $\mathbf{\widehat{U}}_r(r,k)$, and $\omega$. In practice, we only need to consider equations for the first $\frac{1}{2} N +1$ modes, because $\mathbf{U}(r, \theta) \in \mathbb{R}^m$; the $\frac{1}{2} N - 1$ Fourier coefficient vectors $\mathbf{\widehat{U}}(r, k)$ for the modes $k = \frac{1}{2} N + 1, \dots, N - 1$ are the complex conjugates of (non-zero) modes $k = 1, \dots, \frac{1}{2} N-1$. Furthermore, $\mathbf{\widehat{U}}(r,0)$ and $\mathbf{\widehat{U}}(r,\frac{1}{2}N)$ are real. Hence, we are interested in solutions to a system of $2 \, m \, N + 1$ ODEs for $m \, (\frac{1}{2} N + 1)$ real and $m \, (\frac{1}{2} N-1)$ imaginary parts of the Fourier coefficients and their derivatives, as well as the independent variable $r$.

The radial variable $r$ is a time-like variable that fixes a particular
solution segment for this system of ODEs. We consider $r \in [r_0,r_1]$, for
some choices $r_0$ and $r_1$. Together with $\theta \in [0,2\pi)$, this defines
an annulus as the spatial domain for the PDE~\cref{eq:stationary}. The spatial domain resembles a disk if $r_0$ is sufficiently small.
If we were to let 
$r_1\to\infty$, this would amount to a one-dimensional spatial dynamical system connecting 
the core at $r=0$ to a far-field periodic travelling wave in the spirit of \citet{woods1999heteroclinic}. Moreover, computations of so-called boundary sinks, which connect the core to the far field, are useful when investigating spatiotemporal instabilities \citep{Dodson}.

To ensure that system~\cref{eq:autonomousODE} has a well-defined solution for
each of the $\frac{1}{2} N + 1$ Fourier modes, we need to impose $2 \, m \, N + 1$ boundary conditions and an additional phase-pinning condition. As
suggested by \citet{Bordyugov2007}, we impose Neumann (zero radial derivative,
or no flux) boundary conditions on $\mathbf{U}$ on the inner and outer radial
boundary of the chosen annulus in the spatial domain; in terms of the
corresponding Fourier coefficients, these are
 \begin{equation}
  \label{eq:Neumann} 
  \begin{array}{lcr}
    \mathbf{\widehat{U}}^{\prime}(r_0, k) &=& 0, \\
    \mathbf{\widehat{U}}^{\prime}(r_1, k) &=& 0.
  \end{array}
 \end{equation}
To ensure uniqueness of
the computed spiral wave, it is standard to impose a phase-pinning
condition~\citep{barkley1992}. We use the suggestion from~\citep{Bordyugov2007}
and fix the imaginary part of the first Fourier mode (or any other Fourier mode
$k \ge 1$) of one of the variables at $r = r_1$:
 \begin{equation}
  \label{eq:phase}
  {\rm Im}(\widehat{a}(r_1, 1)) = {\rm constant},
 \end{equation}
where we set this constant to $0$. An alternative approach would be to replace this phase-pinning boundary condition with an integral condition. However, the described single boundary condition is qualitatively equivalent and computationally faster. 

\subsection{Starting point for continuation}
\label{sec:StartingSoln}
The solution to the BVP~\cref{eq:autonomousODE}--\cref{eq:phase} is an
$r$-dependent family of Fourier coefficient vectors for the first $\frac{1}{2}N
+ 1$ modes. Using these Fourier coefficients and their complex-conjugates, we
can reconstruct an approximation to the spiral wave defined on the annulus
centred at the origin with radius $r \in [r_0, r_1]$ in the space-time domain.
To obtain such a solution, we follow the continuation approach
from~\citep{Bordyugov2007} that constructs the solution starting from an
(almost) infinitely thin annulus. More precisely, while our goal is to find the
spiral wave on a large disk, that is, with $r_0=0$ and very large $r_1$, we start from a situation where $r_0 \approx r_1$
and allow $r_0$ and $r_1$ to vary, so that $r_0$ becomes very small and,
subsequently, $r_1$ very large. This approach is useful, because the solution
with $r_0 = r_1$ away from $0$ can be approximated by a periodic travelling
wave. Indeed, since the time periodicity of a (rigidly rotating) spiral wave is
equivalent to a rotation in the $(x, y)$-plane, the solution restricted to the
circle with radius $R = r_0 = r_1$ is, in fact, given by a periodic function of
both space and time.

Periodic travelling waves can be found as periodic orbits of the one-dimensional  travelling-frame equation
\begin{equation}
  \label{eq:travellingEqnOriginal}
  \mathbf{f}(\mathbf{U}) + \mathbf{U}_{\xi \xi} - \gamma \, \mathbf{U}_{\xi} = 0,
\end{equation}
which is a stationary solution of the one-dimensional form of PDE~\cref{eq:VectorForm} along the $x$-direction in absence of $y$, formulated with respect to the travelling-frame variable $\xi = x + \gamma \, t$ and the wavespeed $\gamma$ of the periodic travelling wave. Such periodic travelling waves have wavelength $L$ that depends on $\gamma$ \citep{postlethwaite2017}. In \citep{postlethwaite2017, postlethwaite2019,hasan2021spatiotemporal}, we investigated and computed such travelling waves, along with their stability properties, and showed that they, together with the limiting heteroclinic cycles for large wavelength values, provide a mechanism for the existence and stability of spiral waves of system~\cref{eq:PDElaboratory}.

Due to the spatiotemporal symmetry of a (rigidly rotating) spiral wave, a solution to equation~\cref{eq:travellingEqnOriginal} along a radial direction can equivalently be obtained as a solution along the angular direction, formulated in the co-rotating variable $\theta = \phi + \omega \, t$ as
\begin{equation}
  \label{eq:travelling-frame} 
  \mathbf{f}(\mathbf{U}) + \frac{1}{R^2} \mathbf{U}_{\theta \theta} - \omega \mathbf{U}_{\theta} = 0,  
\end{equation}
where we assume that $r = R$ is constant. Since $\theta$ is an angular coordinate, it is natural to impose periodic boundary conditions
 \begin{displaymath}
  \mathbf{U}(2 \pi) = \mathbf{U}(0) \quad  \mbox{and} \quad \mathbf{U}_{\theta}(2 \pi) = \mathbf{U}_{\theta}(0).
 \end{displaymath}
We remark that periodic solutions of  \cref{eq:travelling-frame}  coincide with one-dimensional periodic orbits of \cref{eq:travellingEqnOriginal} with wavespeed $\gamma=R\omega$ and wavelength 
$L=2\pi{R}$.

Translated in terms of Fourier modes, we assume that $\mathbf{\widehat{U}}^{\prime}(r, k) \approx 0$ and $\mathbf{\widehat{U}}^{\prime \prime}(r, k) \approx 0$ in system~\cref{eq:autonomousODE}. In other words, the periodic solution of~\cref{eq:travelling-frame} is an approximate solution of
 \begin{displaymath}
   \mathbf{\widehat{f}}(k) - \frac{k^2}{R^2} \, \mathbf{\widehat{U}}(k) - i \, k \, \omega \, \mathbf{\widehat{U}}(k) = 0,
 \end{displaymath}
defined on a very thin annulus with $r_1 = r_0 + \varepsilon$, for $0 <
\varepsilon \ll 1$.

 \begin{figure}[t!] 
   \centering
   \includegraphics[scale=1.0]{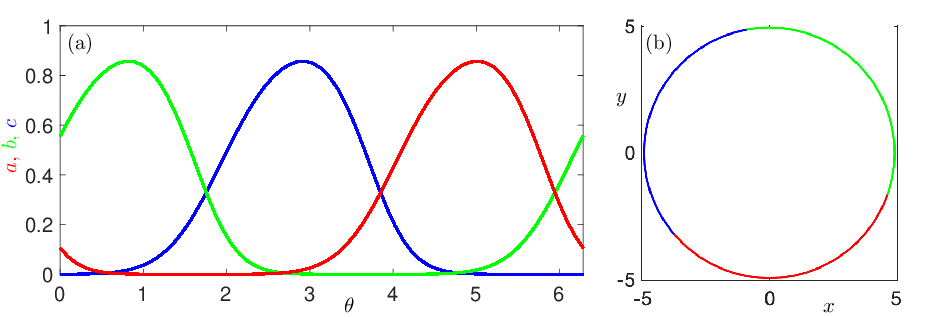}
     \caption{\label{fig:thin}
     The first step for continuing spiral waves of system~\cref{eq:PDElaboratory} with $\sigma = 3.2$, $\zeta = 0.8$, $\omega \approx 0.3346$ and radius $R = 5$. Panel~(a) shows the $\theta$-periodic solution for system~\cref{eq:travelling-frame} and panel~(b) the corresponding distribution of the three populations $a$ (red), $b$ (green) and $c$ (blue) given by system~\cref{eq:autonomousODE} and defined on a thin annulus in the $(x, y)$-plane of width $\varepsilon = 0.002$ centred around $r = R$.}
\end{figure}
%

\Cref{fig:thin} shows the starting solution for the BVP~\cref{eq:autonomousODE}--\cref{eq:phase} with $\sigma=3.2$ and $\zeta=0.8$.  Panel~(a) shows the periodic solution of system~\cref{eq:travelling-frame} with period $2 \pi$ for $R = 5$, which exists when $\omega~\approx~0.3346$. This periodic solution is created in a Hopf bifurcation at $\omega~\approx~0.2640$.
Panel~(b) shows the distribution of the one-dimensional periodic solution along a thin two-dimensional annulus; in Fourier space, this solution is given by approximately constant Fourier coefficient vectors for all $r \in [r_0, r_1]$. Here, we used $\varepsilon = 0.002$ and set
\begin{displaymath}
  r_0 = R - \tfrac{1}{2} \varepsilon = 4.999 \quad \mbox{and} \quad r_1 = R + \tfrac{1}{2} \varepsilon = 5.001.
\end{displaymath}
The red, green and blue regions represent the spatial dominance of $a$, $b$ and $c$, respectively.

The spiral wave can now be found on a larger domain by continuation of this solution to the BVP~\cref{eq:autonomousODE}--\cref{eq:phase}, where we first decrease $r_0$ and then increase $r_1$; the resulting spiral wave for $r \in [0.01, 600]$ is shown in \cref{fig:LargeSW}. For the example system~\cref{eq:PDElaboratory}, where $m = 3$, this means the one-parameter continuation of a system of $2 \, m \, N + 1 = 6 N  + 1$ real ODEs with the same number of unknowns. To obtain a good approximation to the spiral wave, it is imperative that a sufficiently large number $N$ of Fourier modes is computed. 
The appropriate choice for $N$ depends on the system. The almost two-fold reduction of the BVP, obtained by taking into account that real-valued solution vectors must have complex-conjugate pairs of Fourier coefficients, still leaves us with a very large system of equations. Since we wish to explore cases where the far-field dynamics limits onto the one-dimensional TW solution of~\cref{eq:travelling-frame}, and since these solutions are associated with heteroclinic cycles, potentially with sharp transitions, one may need to work with many Fourier modes.  Nevertheless, for the examples considered in this paper, we find that $N = 60$ is sufficient. 

 We note here that \citet*{barkley1992} and \citet*{bar2003} constructed the initial solution by descretising the BVP into a large-scale system of Fourier-space ODEs (with $121 \times 256$ and $51 \times 128$ grid points, respectively) and integrating the ODEs until a spiral wave is obtained. This process is computationally expensive and requires a choice of parameters for which the existence of a spiral wave is known \emph{a priori}. In our study, we obtain the initial guess more reliably and efficiently by finding periodic orbits in a six-dimensional travelling-frame system of ODEs. \citet{Dodson} solve the stationary problem for spiral waves using the MATLAB built-in function \texttt{fsolve} at each continuation step. In contrast, we perform all computations with the BVP set-up supported by the pseudo-arclength continuation software AUTO~\citep{Auto_original,Doedel}; this is also the approach taken by~\cite{Bordyugov2007}.

The next steps in developing the method are to take advantage of the cyclic 
symmetry of this problem, and to address the issue of the boundary condition at 
the inner edge of the annulus, which needs to be changed once the inner
radius~$r_0$ approaches zero.

\subsection{Exploiting cyclic symmetry}
\label{sec:symmetry}

For systems of the form~\cref{eq:VectorForm} satisfying \Cref{assumption:symmetry}, the cyclic symmetry of the nonlinear kinetic term provides an opportunity for further reduction of the number of Fourier coefficients that must be computed in the BVP~\cref{eq:autonomousODE}--\cref{eq:phase} to obtain an $N$-mode approximation of a spiral wave solution. Indeed, the Fourier coefficients also have cyclic symmetry, because the Fourier transformation is an equivariant of this symmetry group. Note that the discretisation of the Fourier transform breaks this equivariance unless $N$ is divisible by $m$. We explain how to implement the symmetry reduction for this restricted class of reaction-diffusion systems with the three-species model~\cref{eq:PDElaboratory} as an example; hence, we assume that $N$ is an integer multiple of $2m=6$.

The cyclic symmetry of system~\cref{eq:PDElaboratory}, defined as the permutation $(a,b,c) \mapsto (b,c,a)$, introduces an additional phase-shift invariance of its spiral waves in the azimuthal direction, namely,
\begin{displaymath}
  \begin{array}{lll}
    b(r, \theta) = a(r, \theta - \frac{2 \pi}{3}), \\
    c(r, \theta) = a(r, \theta + \frac{2 \pi}{3}),
  \end{array}
\end{displaymath}
for all $r$ and $\theta$. 
Hence, it is sufficient to consider only the first component of the co-rotating frame equation~\cref{eq:stationary}, which is given by
\begin{displaymath}
  f_1 \left( a(r, \theta),\, a(r, \theta - \tfrac{2\pi}{3}),\, a(r, \theta + \tfrac{2\pi}{3}) \right) +
  a_{rr}+ \frac{1}{r} \, a_{r} + \frac{1}{r^2} \, a_{\theta \theta} - \omega \, a_{\theta} = 0,
\end{displaymath}
where $f_1$ is the first component of the kinetics function $\mathbf{f}$.
This is related to the approach by \citep{wulff2006} to compute periodic orbits in ODEs with cyclic symmetry.
The co-rotating frame equation in Fourier space for each mode $k$ is then given by
\begin{equation}
  \label{eq:ReducedProblem}
  \widehat{f}_1(r, k) + \widehat{a}_{rr}(r, k) + \frac{1}{r} \, \widehat{a}_{r}(r, k)
                - \frac{k^2}{r^2} \, \widehat{a}(k) - i \, k \, \omega \, \widehat{a}(r, k) = 0.
\end{equation} 
Here similarly, $\widehat{f}_1(r, k)$ is the $k$th coeffcient of the discretised Fourier approximation of $f_1$, which is the first component of the Fourier coefficient vector~\cref{eq:nonlinearTerm}, given by
\begin{displaymath}
  \widehat{a}(r, k) - [\widehat{a} \ast \widehat{a}](r, k) - (1 + \sigma +\zeta) \, [\widehat{a} \ast \widehat{b}](r, k) - (1-\zeta) \,  [\widehat{a} \ast \widehat{c}](r, k).
\end{displaymath}
Since $N$ is assumed to be divisible by $m = 3$, we can express $\widehat{b}(r, k)$ and $\widehat{c}(r, k)$ also in terms of $\widehat{a}(r, k)$, namely,
\begin{displaymath}
  \widehat{b}(r, k) = \widehat{a}(r, k) \, e^{-2\pi \,  i \, k \,  / 3} \quad \mbox{and} \quad
  \widehat{c}(r, k) = \widehat{a}(r, k) \, e^{2\pi \, i \, k \,  / 3}.
\end{displaymath}
As a consequence, we now only need to solve $2 N + 1$ real ODEs.

The boundary conditions for the reduced BVP are just the first components of each of the two conditions in~\cref{eq:Neumann}, with the same phase-pinning condition~\cref{eq:phase}. Hence, they are given by
\begin{equation}
  \label{eq:reducedBCs} 
  \begin{array}{rll}
    \widehat{a}^{\prime}(r_0, k) &=& 0, \\
    \widehat{a}^{\prime}(r_1, k) &=& 0, \\
    \text{Im}(\widehat{a}(r_1, 1)) &=& 0.
\end{array}
\end{equation}
In general, for a system with $m$ species, the number of real differential equations is reduced by almost a factor $m$, from $2 \, m \, N + 1$ to $2 \, N + 1$. In other words, the increase in the number of species $m$ in system~\cref{eq:VectorForm} has no effect on the size of the BVP and only adds a minor extra cost to the computation because of the additional convolutions in the first component of the Fourier coefficient vector $\mathbf{\widehat{f}}$.
We note here also that a similar reduction can be obtained when computing spiral waves that obey a different permutation symmetry, involving a subset of the $m$ species, but the efficiency gains in such instances will be smaller than a factor $m$.

\subsection{Growing the spiral wave solution by continuation}
%
\begin{figure}[t!] 
  \centering
  \includegraphics[scale=1.0]{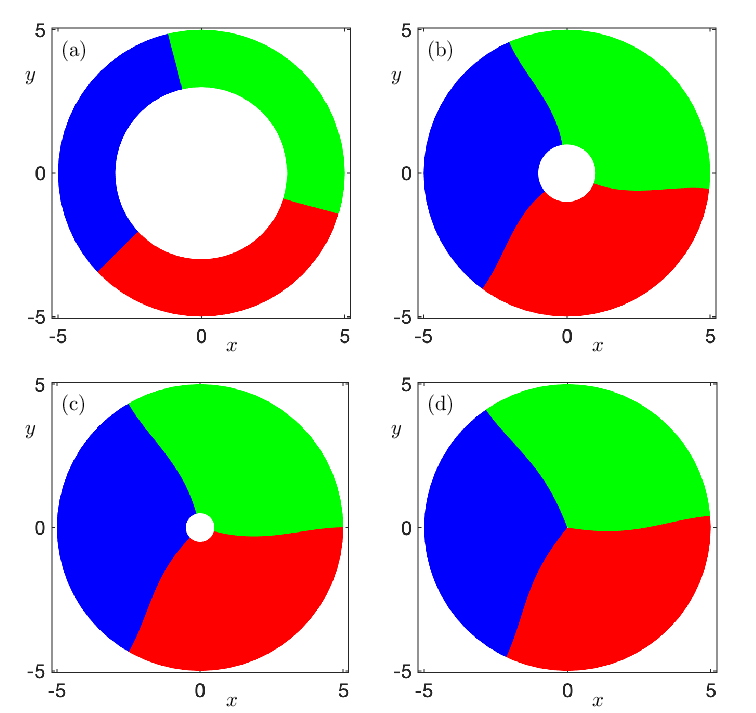}
  \caption{\label{fig:shrinkHole}
    Illustration of the continuation steps taken to compute a spiral wave for system~\cref{eq:PDElaboratory} with $\sigma=3.2$ and $\zeta=0.8$ by varying the inner radius $r_0$. Starting from the solution shown in \cref{fig:thin}(b) with $N = 60$, panels~(a)--(d) show the solutions to the BVP~\cref{eq:ReducedProblem}--\cref{eq:reducedBCs} on the $(x, y)$-plane at a fixed moment in time, where $r_0 = 3$, $r_0 = 1$, $r_0 = 0.5$, and $r_0 = 0.01$, respectively.}
\end{figure}
%
We use the substantially reduced BVP~\cref{eq:ReducedProblem}--\cref{eq:reducedBCs} to compute a spiral wave for system~\cref{eq:PDElaboratory} with $\sigma=3.2$ and $\zeta=0.8$. Recall that we already found a starting solution for these parameter values in the form of a periodic travelling wave interpreted as a solution on a very thin annulus with inner radius $r_0  = 4.999$ and outer radius $r_1 = 5.001$; see \cref{fig:thin}. We set $N = 60$ and initialise the continuation with the Fourier coefficients corresponding to this travelling wave solution. We then continue the BVP~\cref{eq:ReducedProblem}--\cref{eq:reducedBCs} in the direction of decreasing $r_0$ while keeping $r_1$ fixed and treating the angular frequency $\omega$ as a free parameter. \Cref{fig:shrinkHole} shows four snapshots of this continuation process, where we plot the solution in the $(x, y)$-plane with $r_0 = 3$ in panel~(a), $r_0 = 1$ in panel~(b), $r_0 = 0.5$ in panel~(c), and $r_0 = 0.01$ in panel~(d), which is the last step after which we stopped the continuation. Panel~(d) illustrates that the computed solution can be viewed as a computation on a disk of radius $r_1 = 5.001$ in the $(x, y)$-plane, even though the domain is, in reality, an annulus, and the inner boundary condition needs to be considered before $r_0$ can be decreased further. It is also clear that $r_1$ is too small to show the true nature of the spiral wave.
\begin{figure}[t!] 
  \centering
  \includegraphics[scale=1.0]{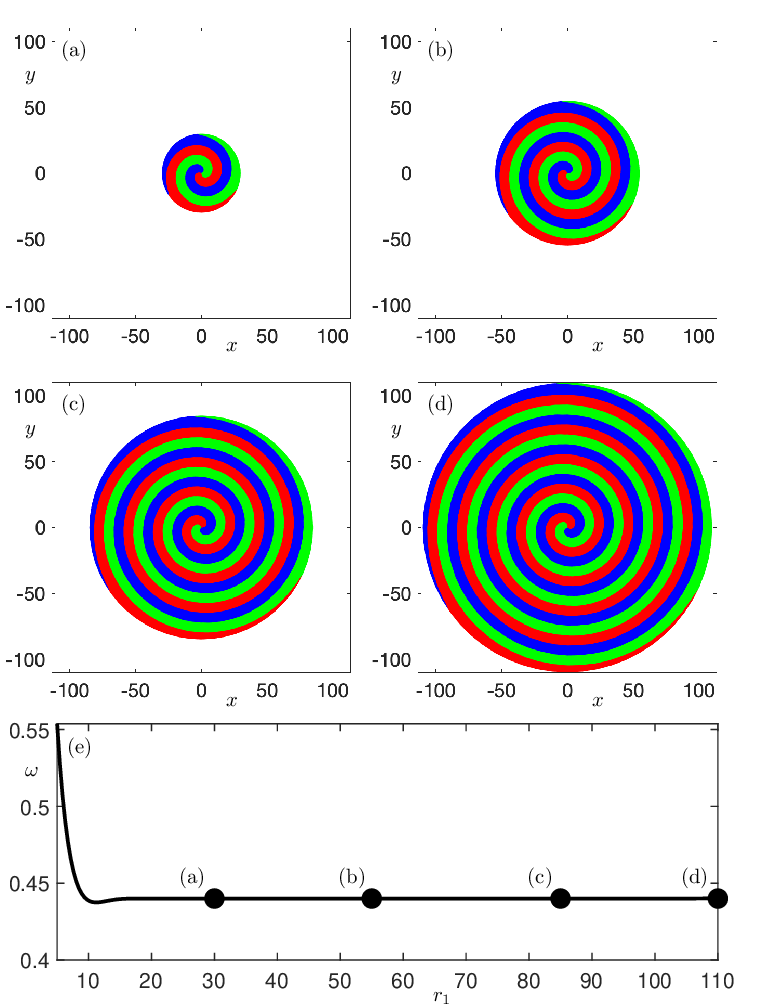}
  \caption{\label{fig:r110}
    Illustration of the continuation steps taken to compute a spiral wave for system~\cref{eq:PDElaboratory} with $\sigma=3.2$, $\zeta=0.8$ and $r_0=0.01$ by varying the outer radius $r_1$. Starting from the solution shown in \cref{fig:shrinkHole}(d) with $N = 60$, panels~(a)--(d) show the solutions to the BVP~\cref{eq:ReducedProblem}--\cref{eq:reducedBCs} on the $(x, y)$-plane at a fixed moment in time, where $r_1= 30$, $r_1= 55$, $r_1=85$, and $r_1=110$, respectively. Panel~(e) shows the angular frequency $\omega$ versus $r_1$.}
\end{figure}

The size of the annular spatial domain can be increased with a second continuation step where we keep $r_0$ fixed and increase  $r_1$ to obtain a spiral wave solution on a sufficiently large domain. \Cref{fig:r110} shows the evolution of the spiral wave as $r_1$ is increased. As in \cref{fig:shrinkHole}, we show each stage of the continuation step by plotting the spiral wave in the $(x, y)$-plane at a fixed moment in time; here, $r_1 = 30$ in panel~(a), $r_1 = 55$ in panel~(b), $r_1 = 85$ in panel~(c), and $r_1 = 110$ in panel~(d). Panel~(e) illustrates how the angular frequency $\omega$ varies during the continuation; here, we plot $\omega$ versus the outer radius $r_1$. Note that the initial value $\omega \approx 0.3346$ for $[r_0, r_1] = [4.999, 5.001]$ has increased to {$\omega \approx 0.5537$ as $r_0$ was decreased to $0.01$. As the outer radius of the spiral is increased from $r_1 = 5.001$, the angular frequency $\omega$ drops again sharply, and then levels off very quickly to $\omega \approx 0.4400$.
We stopped the continuation when $r_1 = 600$; the spiral wave for this value is shown in \cref{fig:LargeSW}(a) and $\omega \approx 0.4400$ here as well.

\subsection{Boundary conditions at the core: spiral waves on a full disk}
\label{sec:BCs}
In the initial thin annulus, we used Neumann boundary conditions in the radial 
direction in order to have an initial solution that is approximately independent 
of~$r$. As the inner radius~$r_0$ is decreased towards zero, these boundary 
conditions are no longer appropriate since they lead to singularities in 
the Laplacian. To see this, consider one Fourier coefficient $\widehat{a}(r,k)$ with wavenumber $k$; for small $r$, this coefficient will behave as
 \begin{displaymath}
 \widehat{a}(r,k) \sim r^j e^{ik\phi},
 \end{displaymath}
where $j\geq0$ is an integer. The corresponding Laplacian is
 \begin{displaymath}
 \nabla^2 r^j e^{ik\phi} = (j^2 - k^2) r^{j-2}  e^{ik\phi}.
 \end{displaymath}
The Laplacian of~$\widehat{a}(r,k)$ should be well behaved at any point on the disk, so at 
$r=0$, we must have either $j^2 - k^2=0$ or $j\geq2$. 
Therefore, for the zeroth
Fourier mode ($k=0$), we have $j=0$ or $j\geq2$, that is, the radial 
derivative of the zeroth mode is zero at $r=0$. For the first
Fourier mode ($k=1$), we can have $j=1$ or $j\geq2$, so the first Fourier mode is zero at $r=0$. Similarly, for $k\geq2$, we must have $j\geq2$.
Imposing Neumann boundary
conditions~\cref{eq:Neumann}, where $\mathbf{\widehat{U}}^{\prime}(r_0, k) = 0$ for all $k$, allows for the possibility that $\widehat{a}(0,k) \neq 0$ for all $k$, which leads to an infinite Laplacian for $k \geq 1$ in the limit $r_0 \to 0$.
Hence, in order to continue to $r_0=0$ without a singularity in the Laplacian, we require a Neumann boundary condition only for $k=0$ and impose Dirichlet boundary conditions at $r = r_0$ for all other
$k$:
 \begin{equation}
  \label{eq:DirichletBCs}
  \begin{array}{rcll}
    \mathbf{\widehat{U}}^{\prime}(r_0, 0) &=& 0, & \\
    \mathbf{\widehat{U}}(r_0, k)              &=& 0, & \mbox{for } k > 0.
  \end{array}
 \end{equation}
As $r_0\to0$, the cyclic symmetry implies that $a$, $b$ and $c$ all take on the same value at $r=0$: we denote this common value by~$\mu_0$.
Note that these boundary conditions translate to a Dirichlet boundary condition in physical space, namely, $\mathbf{U}(0, \theta)= \mu_0  $.
These boundary conditions must be coupled with the phase condition   \cref{eq:phase}.

As was the case for the
Neumann boundary conditions~\cref{eq:Neumann}, the modified Dirichlet boundary conditions~\cref{eq:DirichletBCs} are also readily implemented for the reduced
BVP~\cref{eq:ReducedProblem}--\cref{eq:reducedBCs}, because they translate
directly to conditions for the component $\widehat{a}$ of
$\mathbf{\widehat{U}}$ only.
Rather than set up the modified BVP from scratch (which would not work in the 
initial thin annulus), we continue the spiral
wave that was already computed as a solution to
BVP~\cref{eq:ReducedProblem}--\cref{eq:reducedBCs} with Neumann boundary
conditions to that of a BVP with the modified Dirichlet boundary conditions via
a homotopy step. We define, for $k > 0$, the hybrid boundary conditions
 \begin{displaymath}
  (1 - \lambda) \, \mathbf{\widehat{U}}^{\prime}(r_0, k) + \lambda \,  \mathbf{\widehat{U}}(r_0, k) = 0,
 \end{displaymath}
where we initially set $\lambda = 0$. We then continue the known solution from
$\lambda=0$ to $\lambda=1$, while keeping $\omega$ as a free parameter.
Following this, we can continue $r_0$ to zero resulting in a spiral wave on a full disk with radius $R=r_1$.

\section{Comparison of computational results}
\label{sec:comparethreespecies}

\begin{figure}[t] 
  \centering
  \includegraphics[scale=1.0]{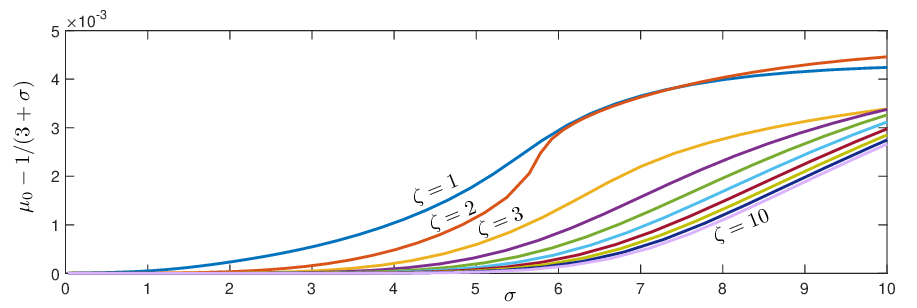}
  \caption{\label{fig:Dirichlet}
    Dependence on the parameter $\sigma$ of the value $a = b = c = \mu_0$ at
    the core of the spiral wave for system~\cref{eq:PDElaboratory} computed
    with radius $R = 30$ and $N=60$ Fourier modes using Dirichlet and Neumann boundary  conditions~\cref{eq:DirichletBCs}. Shown are results from numerical
    continuation with respect to $\sigma$ for fixed $\zeta = 1, 2, \dots 10$.
    The curves illustrate the difference $\mu_0 - \frac{1}{3+\sigma}$ for each of
    the ten different values of $\zeta$.}
 \end{figure}

To test the accuracy and efficiency of our modified method, we compare the
results obtained from the reduced
BVP~\cref{eq:ReducedProblem}--\cref{eq:reducedBCs} with direct simulations of
the laboratory-frame PDEs~\cref{eq:PDElaboratory}. We note that the full
BVP~\cref{eq:autonomousODE}--\cref{eq:phase} produces solutions with the same cyclic
symmetry as assumed in the reduced
BVP~\cref{eq:ReducedProblem}--\cref{eq:reducedBCs}, so the results are
identical, with a substantial gain in efficiency.

For the BVP method, we compute spiral waves on a disk with radius $R=30$ and continue
them with respect to $\omega$, $\sigma$ and~$\zeta$. The PDE simulations were conducted
on a large square box in the $(x, y)$-plane of size $2000 \times 2000$;
see~\citep{postlethwaite2017,hasan2021spatiotemporal} for precise details on this
computational set-up.

\subsection{Values of $a$, $b$ and $c$ at the core}

We have defined $\mu_0$ to be the common value of $a$, $b$ and $c$ at the 
core $r=0$. 
Results from direct PDE simulation show that this common value is close to (but not equal to) the 
coexistence equilibrium point value, 
$a=b=c=\frac{1}{3+\sigma}$~\citep{Szczesny2014,postlethwaite2017}.
 Our continuation scheme gives the same result for the values of $a$, $b$ and $c$ at the core as the direct simulation of the PDEs.
In \cref{fig:Dirichlet} we compare $\mu_0$ to $\frac{1}{3+\sigma}$. 
The continuation results reveal that
$\mu_0\to\frac{1}{3+\sigma}$ as $\sigma\to0$, and more quickly for smaller~$\zeta$. Even for larger~$\sigma$, the
difference is only of order~$10^{-3}$.
We note that this comparison is not possible without the final step of changing the boundary conditions and continuing to $r_0=0$.

\subsection{Angular frequency of spiral waves}
\label{sec:angular}

\begin{figure}[t!] 
  \centering
  \includegraphics[scale=1.0]{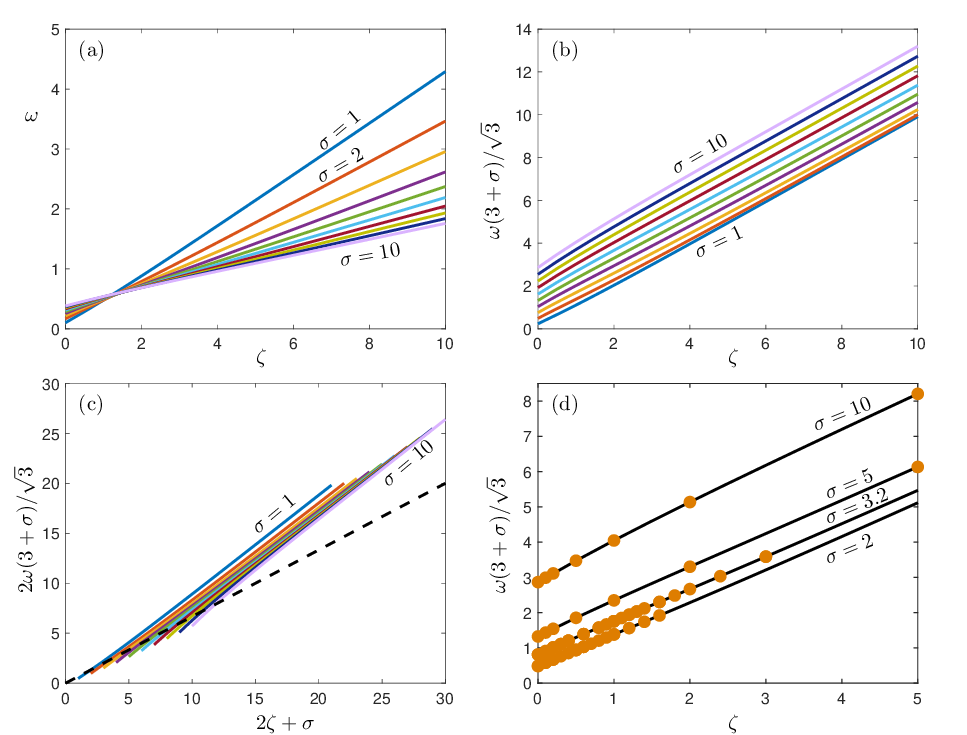}
  \caption{\label{fig:comparison}
    Dependence of the angular frequency $\omega$ on parameters $\zeta$ and $\sigma$ for spiral waves of system~\cref{eq:PDElaboratory} computed on a disk with radius $R = 30$ and $N=60$ Fourier modes.
  Panels (a)$-$(c) show the curves for fixed $\sigma = 1, 2, \dots, 10$ in three different projections.
Panel (d) shows a comparison between the continuation results (black curves) and direct simulations of the laboratory-frame equations~\cref{eq:PDElaboratory} (orange dots) for $\sigma = 2,3.2,5.0$ and $10.0$.
}
\end{figure}
%

\cite{postlethwaite2017} conjectured that the 
angular frequency of the spiral waves was related to the imaginary part of the 
complex eigenvalue at the coexistence equilibrium point. In particular, they 
suggested a linear relationship
 \begin{equation}
 \label{eq:PR2017}
 \omega \sim \frac{2}{3} \times \frac{\sqrt{3}(\sigma+2\zeta)}{2(3+\sigma)},
 \end{equation}
where the $\frac{2}{3}$ pre-factor was obtained by fitting a straight line
to PDE simulations in the range $0.1\leq\sigma\leq20$ and $0\leq\zeta\leq10$. 
The spiral angular frequency in the simulations was only estimated at parameter values 
where reasonably large stable spirals could be found.

To test this conjecture, we computed the angular frequency~$\omega$ by 
continuation. 
\Cref{fig:comparison} shows the relation between the angular frequency $\omega$ and system parameters $\sigma$ and $\zeta$.
Panel~(a) shows this relation in the $(\zeta,\omega)$-plane for $\sigma=1,2, \dots, 10$. 
By inspection, we find that all curves are almost linear with an approximate slope of $\sqrt{3} \,
(3+\sigma)^{-1}$ rather than the conjectured $\frac{2}{3} \sqrt{3} \,
(3+\sigma)^{-1}$ as postulated by \citet{postlethwaite2017}.
This is more evident in panel (b) where we plot the same curves with $\omega$ rescaled by a factor of $\sqrt{3} \,
(3+\sigma)^{-1}$; note that the slopes are all about 1 in this projection.

\Cref{fig:comparison}(c) shows the same curves with $2\omega(3+\sigma)/\sqrt{3}$ plotted 
against $2\zeta+\sigma$.
Our continuation results, over a full range of $\sigma$ and $\zeta$ (each between 1 and 10)
contradict the conjecture of \citet{postlethwaite2017}.
The observations
of \citet{postlethwaite2017} were predominantly based on relatively small values
of $\sigma$ and $\zeta$, where large enough stable spirals could be found and 
their frequencies measured.
There is a rough agreement for small $\sigma$ and $\zeta$, but the discrepancies for larger values are quite
significant; see~\citep[fig.~4]{postlethwaite2017}.
\Cref{fig:comparison}(d) demonstrates that the continuation results (black curves) and direct simulations (orange dots) agree extremely well.

We find that the angular frequency is better approximated by the linear function $$\omega
\approx \frac{\sqrt{3} \, (\frac{3}{2} \sigma + 2 \zeta)}{2(3+\sigma)};$$ 
compare with  \cref{eq:PR2017}. While this value is still similar to the imaginary part of
the coexistence equilibrium, the difference is not just a factor
of~$\frac{2}{3}$. On the other hand, our updated value still indicates linear dependency on $\zeta$ and the slope $\sqrt{3} \, (3+\sigma)^{-1}$ appears to be valid over a large range of values for both $\sigma$ and~$\zeta$.

\section{Five-species model}
\label{sec:5species}
In this section, we illustrate the gains in computational efficiency obtained when taking into account the cyclic symmetry of spiral waves in larger heteroclinic networks.
Heteroclinic networks of five competing species in an ODE form was first studied by \cite{field1992} and then further investigated in \cite{podvigina2013}, \cite{afraimovich2016} and \cite{bayliss2020}.
Just as for the three-strategy Rock--Paper--Scissors game, the extended game of Rock--Paper--Scissors--Lizard--Spock \citep{Kass1995} can be viewed as a population model for competing species; each population density represents the probability density of winning the game when consistently playing the same strategy, e.g., always Rock, or always Scissors. Hence, the game can be modelled as a PDE of the form~\cref{eq:VectorForm} using similar equations as in system~\cref{eq:PDElaboratory}, but with $m = 5$ species. More precisely, we consider 
\begin{equation}
  \label{eq:5speciesPDE} 
  \left\{ \begin{array}{rrlll}
	    \dot{a} &=&  a \, (1 - \rho - (\sigma + \zeta) \, (b + p) + \zeta \, (c+q)) + \nabla^2 a,	 \\
	    \dot{b} &=&  b \, (1 - \rho - (\sigma + \zeta) \, (c + q) + \zeta \, (p+a)) + \nabla^2 b,	 \\
            \dot{c} &=&  c \, (1 - \rho - (\sigma + \zeta) \, (p + a) + \zeta \, (q+b)) + \nabla^2 c,      \\
	    \dot{p} &=&  p \, (1 - \rho - (\sigma + \zeta) \, (q + b) + \zeta \, (a+c)) + \nabla^2 p,	 \\
            \dot{q} &=&  q \, (1 - \rho - (\sigma + \zeta) \, (a + c) + \zeta \, (b+p)) + \nabla^2 q,
  \end{array} \right.
\end{equation}
where $p(t,x,y)$ and $q(t,x,y)$ are two additional (non-dimensionalised)
species and $\rho = a + b + c + p + q$. Here, we assume that the removal rate
$\sigma$ and replacement rate $\zeta$ are the same for all species interactions.
Spatiotemporal patterns for five-species reaction-diffusion systems have been obtained via stochastic simulations and studied in~\citep{hawick2011,cheng2014,
 park2017},
 but a precise continuation or computation of the spiral wave has not been attempted before.

 }

\begin{table}[t!] 
  \centering
  \begin{tabular}{|c|c|c|c|c|c|}
    \hline
    \bf{Equilibrium}  & $\xi_1$     & $\xi_2$     & $\xi_3$     & $\xi_4$     & $\xi_5$     \\ \hline
    $(a,b,c,p,q)$ & $(1,0,0,0,0)$ & $(0,1,0,0,0)$ & $(0,0,1,0,0)$ & $(0,0,0,1,0)$ & $(0,0,0,0,1)$ \\ \hline
    \bf{Equilibrium}   & $\eta_1$     & $\eta_2$     & $\eta_3$     & $\eta_4$     & $\eta_5$     \\ \hline
    $(a,b,c,p,q)$ & $(s,s,0,0,s)$ & $(s,s,s,0,0)$ & $(0,s,s,s,0)$ & $(0,0,s,s,s)$ & $(s,0,0,s,s)$ \\ \hline
  \end{tabular}
  \caption{\label{tb:equilibria}
    Equilibrium points in the heteroclinic network of system~\cref{eq:5speciesPDE}, where $s = (3+\sigma)^{-1}$.}
\end{table}
%
System~\cref{eq:5speciesPDE} has five equilibria at which only one population ($a$, $b$, $c$, $p$, and $q$, respectively) survives, which we denote by $\xi_i$, where the index $i = 1, 2, 3, 4, 5$ represents the component in the vector $(a, b, c, p, q)$ of the surviving species.
There are also five equilibria $\eta_1-\eta_5$, each of which has three surviving populations; more precisely, at each equilibrium $\eta_i$, components with indices $i-1$, $i$,$i+1 \, ({\rm mod} \, 5)$ survive.
The coordinates of the equilibrium points $\xi_1-\xi_5$ and $\eta_1-\eta_5$ are given in~\cref{tb:equilibria}.
 Note that each triple $\xi_{i-1}$, $\xi_i$, $\xi_{i+1}$ lies in an invariant subspace that contains the equilibrium $\eta_i$. 
Moreover, there exists a coexistence equilibrium with $a=b=c=p=q=1/(5+2\sigma)$, where all five species coexist.
\begin{figure}[t!] 
  \centering
  \includegraphics[scale=1.0]{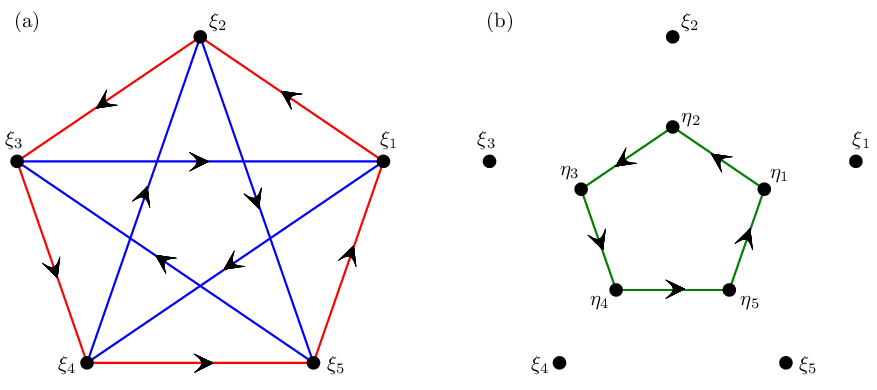}
   \caption{\label{fig:Network}
    Heteroclinic network of system~\cref{eq:5speciesPDE} with nodes $\xi_1,\xi_2,\xi_3,\xi_4$ and $\xi_5$ that represent equilibria at which only one population survives, and $\eta_1,\eta_2,\eta_3,\eta_4$ and $\eta_5$ that represent equilibria at which three populations survive. Panel~(a) highlights heteroclinic cycles $\Gamma_1$ (red) and $\Gamma_2$ (blue) involving all five equilibria $\xi_i$, and panel~(b) shows shows the heteroclinic cycle $\Gamma_3$ (green) involving all five equilibria $\eta_i$, with $i = 1, \dots, 5$; see~\cref{tb:equilibria} for the definition of the equilibria.}
 \end{figure}
%

\Cref{fig:Network} shows all possible five-component heteroclinic cycles of system~\cref{eq:5speciesPDE} in the form of a directed graph, where each node represents an equilibrium and the edges the connecting orbits. As for the three-species model, there exists a heteroclinic cycle that consecutively connects all five equilibria $\xi_i$ with $i = 1, \dots, 5$; this cycle is highlighted in red in \cref{fig:Network}(a) and we refer to it as $\Gamma_1$. There exists a second heteroclinic cycle between these five equilibria: instead of connecting $\xi_i$ to $\xi_{i+1}$ for $i = 1, \dots, 5 \, ({\rm mod} \, 5)$, the connection is from $\xi_i$ to $\xi_{i+3}$ for $i = 1, \dots, 5 \, ({\rm mod} \, 5)$; this cycle is highlighted in blue in \cref{fig:Network}(a) and we refer to it as $\Gamma_2$. There also exists a heteroclinic cycle between the other five equilibria, connecting $\eta_i$ to $\eta_{i+1}$ for $i = 1, \dots, 5 \, ({\rm mod} \, 5)$; this cycle is highlighted in green in \cref{fig:Network}(b) and we refer to it as $\Gamma_3$.

There are many other heteroclinic cycles, some longer and some involving fewer than five equilibria. For example, there exist five heteroclinic cycles in system~\cref{eq:5speciesPDE} between three of the five equilibria $\xi_i$ 
with $i = 1, \dots, 5$. More precisely, these are formed by triples $\xi_{i-1}$, $\xi_i$, $\xi_{i+1}$ that surround an associated coexistence equilibrium $\eta_i$, with $i = 1, \dots, 5 \, ({\rm mod} \, 5)$. Hence, the dynamics is entirely restricted to the lower-dimensional invariant subspace involving these equilibria, which is perfectly described by system~\cref{eq:PDElaboratory}. Here, we focus only on the dynamics of the heteroclinic cycles between five equilibria.

With our reduced BVP method, we take advantage of the cyclic symmetry in the system; in complete analogy to the steps outlined in~\cref{sec:symmetry}, we can express each of the five species in terms of appropriate phase shifts of the first. 
As a starting solution, we determine the associated periodic travelling wave given as a periodic orbit of the one-dimensional travelling-frame equation 
 \begin{equation}
  \label{eq:travellingEqn5species}
  \mathbf{f}(\mathbf{U}) + \mathbf{U}_{\xi \xi} 
  - \gamma \mathbf{U}_{\xi} = 0,
 \end{equation} 
where the solution vector $\mathbf{U}(\xi) = [a(\xi),\, b(\xi),\, c(\xi),\, p(\xi),\, q(\xi)]^T \in \mathbb{R}^m$ now represents $m = 5$ species, the vector $\mathbf{f}$ represents the kinetics terms in system \cref{eq:5speciesPDE}, and $\xi=x+\gamma t$ is the travelling-frame variable.
In the travelling-frame coordinate, we find periodic orbits that are associated to only two of the three heteroclinic cycles between five equilibria, namely, $\Gamma_2$ and $\Gamma_3$, that is, the blue and green cycles highlighted in \cref{fig:Network}. Using these two periodic orbits, we set up two different versions of the reduced BVP~\cref{eq:ReducedProblem}--\cref{eq:reducedBCs}: for cycle $\Gamma_2$, we can exploit the permutation symmetry
\begin{displaymath}
 (a, b, c, p, q) \mapsto (p, q, a, b, c),
\end{displaymath}
and for cycle $\Gamma_3$, we use
\begin{displaymath}
  (a, b, c, p, q) \mapsto (b, c, p, q, a).
\end{displaymath} 
Therefore, the first component of the co-rotating frame equation~\cref{eq:stationary} for cycle $\Gamma_2$ is given by
\begin{eqnarray*}
  f_1 \left( a(r, \theta),\, a(r, \theta - \tfrac{4\pi}{5}),\, a(r, \theta + \tfrac{2\pi}{5}),\, a(r, \theta - \tfrac{2\pi}{5}),\, a(r, \theta + \tfrac{4\pi}{5}) ) \right) + & & \\
  a_{rr}+ \frac{1}{r} \, a_{r} + \frac{1}{r^2} \, a_{\theta \theta} - \omega \, a_{\theta} &=& 0,
\end{eqnarray*}
where
\begin{displaymath}
  f_1(a, b, c, p, q) = a \left(1 - \rho - (\sigma + \zeta) \, (b + p) + \zeta \, (c + q) \right).
\end{displaymath}
Similarly, for cycle $\Gamma_3$, we use
\begin{eqnarray*}
  f_1 \left( a(r, \theta),\, a(r, \theta - \tfrac{2\pi}{5}),\, a(r, \theta - \tfrac{4\pi}{5}),\, a(r, \theta + \tfrac{4\pi}{5}),\, a(r, \theta + \tfrac{2\pi}{5}) ) \right) + & & \\
  a_{rr}+ \frac{1}{r} \, a_{r} + \frac{1}{r^2} \, a_{\theta \theta} - \omega \, a_{\theta} &=& 0.
\end{eqnarray*}
In this way, we reduce the number of real differential equations from $2 \, m \, N + 1 = 10 \, N + 1$ to $2 \, N + 1$,
which is about a five-fold reduction.

The numerical continuation is performed in Fourier space, using equation~\cref{eq:ReducedProblem} for each of the first $\frac{1}{2}N+1$ Fourier modes. The $k$th coefficient of the discretised Fourier approximation of $f_1$ is now given by
\begin{equation}
  \label{eq:nonlinear5species}
\widehat{f}_1(r,k)= \  \widehat{a} - \widehat{a} \ast \widehat{a} - (1 + \sigma + \zeta) \, \widehat{a} \ast (\widehat{b} + \widehat{p}) - (1 - \zeta) \, \widehat{a} \ast (\widehat{c} + \widehat{q}),
\end{equation}
where we dropped the dependence on $(r, k)$ for notational convenience.

The modifications of the set-up for the BVP~\cref{eq:ReducedProblem}--\cref{eq:reducedBCs} are all that is needed to compute spiral waves of system~\cref{eq:5speciesPDE}. Note that it hardly matters that there are now five rather than three species, because we still only compute the evolution of the population density of one of these species; the number of differential equations and boundary conditions are the same for both the three-species and five-species models. 
The only difference in computational cost is the evaluation of five instead of three  nonlinear terms.
\begin{figure}[t!] 
  \centering
  \includegraphics[scale=1.0]{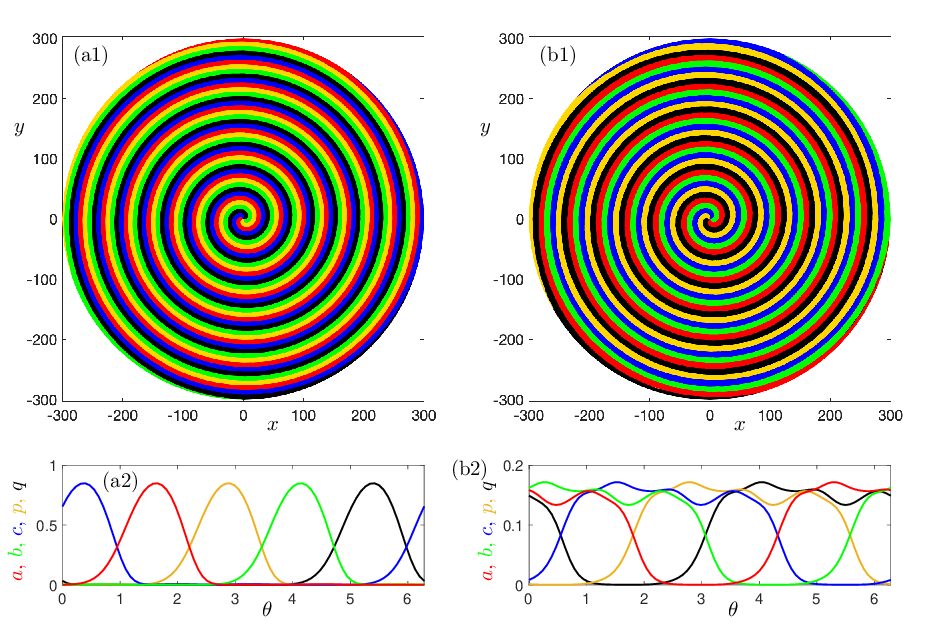}
  \caption{\label{fig:5speciesSW}
    Spiral waves of the five-species system~\cref{eq:5speciesPDE} with $\sigma=3.2$, $\zeta=0.8$ computed on a full disk of radius $R=300$ with $N=80$ Fourier modes. Panels~(a1) and~(b1) show the spiral waves that correspond to cycles $\Gamma_2$ and $\Gamma_3$ with angular velocities $\omega \approx 0.260376$ and $\omega \approx 0.093288$, respectively; panels~(a2) and~(b2) show the corresponding distribution of the five species along outer boundary of the domain.}
\end{figure}

We compute the spiral waves of system~\cref{eq:5speciesPDE} that correspond to the two cycles $\Gamma_2$ and $\Gamma_3$ on a full disk in the $(x, y)$-plane, centred at the origin, with inner radius $r_0 = 0$ and outer radius $r_1 = 300$, by solving the BVP\cref{eq:ReducedProblem}--\cref{eq:reducedBCs} with Dirichlet boundary conditions \cref{eq:DirichletBCs}. The result is shown in \cref{fig:5speciesSW}, where the colours red, green, blue, yellow, and black represent the populations $a$, $b$, $c$, $p$, and $q$, respectively. In the top row, we show the spiral wave in the ($x,y$)-plane with colours indicating the species that is dominant in that region. The bottom row shows the spatial profile along the outermost circle of the computed domain (at $r = 300$); here the horizonal axis is the $2 \pi$-periodic co-rotating variable $\theta$. Panels~(a1) and~(a2) are for the spiral wave associated with the cycle $\Gamma_2$ and panels~(b1) and~(b2) for the spiral wave associated with $\Gamma_3$. Note that the profile shown in panel~(a2) lies close to the heteroclinic cycle $\Gamma_2$. Similarly, the profile shown in panel~(b2) is reminiscent of the heteroclinic cycle $\Gamma_3$; the small oscillations around $\frac{1}{3+\sigma}$ are formed due to the saddle-focus nature of the equilibria $\eta_1 - \eta_5$.
The core of both spiral waves shown in~\cref{fig:5speciesSW} admits a common value for all five species.
This value is close to the coexistence equilibrium  $a=b=c=p=q=1/(5+2\sigma)$, for both spirals, within $10^{-3}$.

\section{Discussion}
\label{sec:discussion}

This paper presents a numerical study of spiral waves in systems of symmetric
heteroclinic networks. As an adaptation of the method introduced by
\cite{Bordyugov2007}, we combined methods from discrete Fourier transforms,
boundary-value problem formulations and group symmetries to generate a
fast algorithm for computing and continuing spiral waves in cyclic dominance
models. The proposed continuation method is an efficient way to explore the
dynamics of spiral waves in symmetric systems. In particular, this method is an
efficient way to explore spiral waves in a large heteroclinic networks of five
(or more) competing species.
 We emphasise that the numerical method presented here can be  modified to compute spiral waves for any reaction-diffusion system. More precisely, one needs to follow the steps presented in \cref{sec:background} but ignore the symmetry exploitation presented in \cref{sec:symmetry}.

We found that the angular frequency of spiral waves in system
\cref{eq:PDElaboratory} was not directly related to the linear frequency of
the equilibrium point at the core, as postulated by \citet{postlethwaite2017}.
The results from the modified continuation method also accurately match results from direct integration. 
We introduced Dirichlet boundary conditions \cref{eq:DirichletBCs} in order to obtain a bounded (non-singular) Laplacian term at the core as $r \to 0$. 
Applying these boundary conditions allow us to compute spiral waves on a full disk instead of an annulus.
What is more, a challenging task in direct simulations is to locate and determine the common value of the variables at the core.  
Using our boundary conditions, we are able to determine the precise location and value at the core and investigate its dynamics.

Spatiotemporal instabilities of spiral waves come in different shapes and forms and may emerge from the far field, the core, or the boundary conditions \citep{Dodson,sandstede2000absolute,SandScheel2020}. 
The Dirichlet boundary conditions \cref{eq:DirichletBCs} proposed in this paper is an essential starting point to compute and identify core instabilities. 
Just as we have done for travelling waves \citep{hasan2021spatiotemporal}, we plan to compute spectra of heteroclinic-induced spiral waves.
 Such computations require a coupling of the stationary problem and the augmented eigenvalue problem. Hence, obtaining linear spectra of spiral waves in large domains is a computationally expensive task.
Our symmetry-based reduction can be applied to reduce the number of ODEs for the augmented eigenvalue problem and thus facilitate fast and efficient computation of spectra of spiral waves in spatiotemporal systems of large heteroclinic networks.
One shortcoming of this approach is that the the augmented eigenvalue problem would have full rank, which makes the stability analysis computationally expensive.
However, one can perform iterative steps to compute the largest eigenvalue or even a selected number of dominant eigenvalues efficiently~\citep{dijkstra2014numerical,saad2011numerical}.



For the five-species system, we were only able to find spiral waves associated with two of the three heteroclinic cycles between five equilibria. At first glance, it may seem like the heteroclinic cycles $\Gamma_1$ and $\Gamma_2$ are of a similar type, but the heteroclinic connections in $\Gamma_1$ are two dimensional, and the heteroclinic connections in $\Gamma_2$ are one dimensional--more precisely, there exists a two-dimensional manifold of connecting orbits between $\xi_1$ and $\xi_2$. It may be that this structural difference has a role to play in the bifurcation of long-period periodic orbits from these heteroclinic cycles.
We plan to investigate this in future work by examining system~\cref{eq:5speciesPDE} with one spatial dimension in the travelling frame coordinates and investigating the heteroclinic bifurcations which occur in the resulting 10-dimensional system of ODEs.


\section*{Acknowledgment}

We are grateful for discussions with Andrus Giraldo, Chris Marcotte, Mauro Mobilia and Stephanie Dodson. 

\section*{Funding}
University of Auckland (FRDF-3714924 to C.R.H., H.M.O. \& C.M.P.); Royal Society Te Ap\=arangi, New Zealand (17-UOA-096 to C.M.P.); London Mathematical Laboratory; Leverhulme Trust (RF-2018-449/9 to A.M.R.).

\bibliographystyle{imamat}
\bibliography{sample}

\begin{thebibliography}{}

\bibitem[Afraimovich et~al., 2016]{afraimovich2016}
Afraimovich, V.~S., Moses, G. {\&} Young, T. (2016)  Two-dimensional
  heteroclinic attractor in the generalized Lotka--Volterra system. {\em
  Nonlinearity}, \textbf{29}(5), 1645.

\bibitem[B{\"a}r et~al., 2003]{bar2003}
B{\"a}r, M., Bangia, A.~K. {\&} Kevrekidis, I.~G. (2003)  Bifurcation and
  stability analysis of rotating chemical spirals in circular domains:
  Boundary-induced meandering and stabilization. {\em Physical Review E},
  \textbf{67}(5), 056126.

\bibitem[Barkley, 1992]{barkley1992}
Barkley, D. (1992)  Linear stability analysis of rotating spiral waves in
  excitable media. {\em Physical Review Letters}, \textbf{68}(13), 2090.

\bibitem[Bayliss et~al., 2020]{bayliss2020}
Bayliss, A., Nepomnyashchy, A. {\&} Volpert, V. (2020)  Beyond
  rock-paper-scissors systems-Deterministic models of cyclic ecological systems
  with more than three species. {\em Physica D: Nonlinear Phenomena},
  \textbf{411}, 132585.

\bibitem[Bordyugov \& Engel, 2007]{Bordyugov2007}
Bordyugov, G. {\&} Engel, H. (2007)  Continuation of spiral waves. {\em Physica
  D: Nonlinear Phenomena}, \textbf{228}(1), 49--58.

\bibitem[Cheng et~al., 2014]{cheng2014}
Cheng, H., Yao, N., Huang, Z.-G., Park, J., Do, Y. {\&} Lai, Y.-C. (2014)
  Mesoscopic interactions and species coexistence in evolutionary game dynamics
  of cyclic competitions. {\em Scientific reports}, \textbf{4}(1), 1--7.

\bibitem[Dijkstra et~al., 2014]{dijkstra2014numerical}
Dijkstra, H.~A., Wubs, F.~W., Cliffe, A.~K., Doedel, E., Dragomirescu, I.~F.,
  Eckhardt, B., Gelfgat, A.~Y., Hazel, A.~L., Lucarini, V., Salinger, A.~G.,
  Phipps, E.~T., Sanchez-Umbria, J., Schuttelaars, H., Tuckerman, L.~S. {\&}
  Thiele, U. (2014)  Numerical bifurcation methods and their application to
  fluid dynamics: analysis beyond simulation. {\em Communications in
  Computational Physics}, \textbf{15}(1), 1--45.

\bibitem[Dodson \& Sandstede, 2019]{Dodson}
Dodson, S. {\&} Sandstede, B. (2019)  Determining the source of period-doubling
  instabilities in spiral waves. {\em SIAM Journal on Applied Dynamical
  Systems}, \textbf{18}(4), 2202--2226.

\bibitem[Doedel, 1981]{Auto_original}
Doedel, E.~J. (1981)  AUTO: A program for the automatic bifurcation analysis of
  autonomous systems. {\em Congr. Numer}, \textbf{30}(3), 265--284.

\bibitem[Doedel et~al., 2007]{Doedel}
Doedel, E.~J., Champneys, A.~R., Dercole, F., Fairgrieve, T.~F., Kuznetsov,
  Y.~A., Oldeman, B., Paffenroth, R., Sandstede, B., Wang, X. {\&} Zhang, C.
  (2007) {\em AUTO-07P: Continuation and bifurcation software for ordinary
  differential equations}.

\bibitem[Field \& Richardson, 1992]{field1992}
Field, M. {\&} Richardson, R. (1992)  Symmetry breaking and branching patterns
  in equivariant bifurcation theory II. {\em Archive for rational mechanics and
  analysis}, \textbf{120}(2), 147--190.

\bibitem[Field \& Swift, 1991]{field1991}
Field, M. {\&} Swift, J.~W. (1991)  Stationary bifurcation to limit cycles and
  heteroclinic cycles. {\em Nonlinearity}, \textbf{4}(4), 1001.

\bibitem[Frey, 2010]{Frey2010}
Frey, E. (2010)  Evolutionary game theory: Theoretical concepts and
  applications to microbial communities. {\em Physica A: Statistical Mechanics
  and its Applications}, \textbf{389}(20), 4265--4298.

\bibitem[Guckenheimer \& Holmes, 1988]{guckenheimer1988}
Guckenheimer, J. {\&} Holmes, P. (1988)  Structurally stable heteroclinic
  cycles. In {\em Mathematical Proceedings of the Cambridge Philosophical
  Society}, volume 103, pages 189--192. Cambridge University Press.

\bibitem[Hasan et~al., 2021]{hasan2021spatiotemporal}
Hasan, C.~R., Osinga, H.~M., Postlethwaite, C.~M. {\&} Rucklidge, A.~M. (2021)
  Spatiotemporal stability of periodic travelling waves in a heteroclinic-cycle
  model. {\em Nonlinearity}, \textbf{34}(8), 5576.

\bibitem[Hawick, 2011]{hawick2011}
Hawick, K. (2011)  Cycles, diversity and competition in
  rock-paper-scissors-lizard-spock spatial game agent simulations. In {\em
  Proceedings on the International Conference on Artificial Intelligence
  (ICAI)}, page~1. The Steering Committee of The World Congress in Computer
  Science, Computer.

\bibitem[Jackson \& Buss, 1975]{Jackson1975}
Jackson, J. {\&} Buss, L. (1975)  Alleopathy and spatial competition among
  coral reef invertebrates. {\em Proceedings of the National Academy of
  Sciences}, \textbf{72}(12), 5160--5163.

\bibitem[Kass \& Bryla, 1995]{Kass1995}
Kass, S. {\&} Bryla, K. (1995) {\em Rock paper scissors Spock lizard}.

\bibitem[Kerr et~al., 2002]{Kerr2002}
Kerr, B., Riley, M.~A., Feldman, M.~W. {\&} Bohannan, B.~J. (2002)  Local
  dispersal promotes biodiversity in a real-life game of rock--paper--scissors.
  {\em Nature}, \textbf{418}(6894), 171--174.

\bibitem[Kirkup \& Riley, 2004]{Kirkup2004}
Kirkup, B.~C. {\&} Riley, M.~A. (2004)  Antibiotic-mediated antagonism leads to
  a bacterial game of rock--paper--scissors in vivo. {\em Nature},
  \textbf{428}(6981), 412--414.

\bibitem[Krupa, 1997]{krupa1997}
Krupa, M. (1997)  Robust heteroclinic cycles. {\em Journal of Nonlinear
  Science}, \textbf{7}(2), 129--176.

\bibitem[May \& Leonard, 1975]{May1975}
May, R.~M. {\&} Leonard, W.~J. (1975)  Nonlinear aspects of competition between
  three species. {\em SIAM journal on applied mathematics}, \textbf{29}(2),
  243--253.

\bibitem[Park et~al., 2017]{park2017}
Park, J., Do, Y., Jang, B. {\&} Lai, Y.-C. (2017)  Emergence of unusual
  coexistence states in cyclic game systems. {\em Scientific reports},
  \textbf{7}(1), 1--11.

\bibitem[Podvigina, 2013]{podvigina2013}
Podvigina, O. (2013)  Classification and stability of simple homoclinic cycles
  in $\mathbb{R}^5$. {\em Nonlinearity}, \textbf{26}(5), 1501.

\bibitem[Postlethwaite \& Rucklidge, 2017]{postlethwaite2017}
Postlethwaite, C. {\&} Rucklidge, A. (2017)  Spirals and heteroclinic cycles in
  a spatially extended rock-paper-scissors model of cyclic dominance. {\em EPL
  (Europhysics Letters)}, \textbf{117}(4), 48006.

\bibitem[Postlethwaite \& Rucklidge, 2019]{postlethwaite2019}
Postlethwaite, C.~M. {\&} Rucklidge, A.~M. (2019)  A trio of heteroclinic
  bifurcations arising from a model of spatially-extended
  Rock--Paper--Scissors. {\em Nonlinearity}, \textbf{32}(4), 1375.

\bibitem[Proctor \& Jones, 1988]{proctor1988}
Proctor, M.~R. {\&} Jones, C.~A. (1988)  The interaction of two spatially
  resonant patterns in thermal convection. Part 1. Exact 1: 2 resonance. {\em
  Journal of Fluid Mechanics}, \textbf{188}, 301--335.

\bibitem[Reichenbach et~al., 2007]{reichenbach2007}
Reichenbach, T., Mobilia, M. {\&} Frey, E. (2007)  Mobility promotes and
  jeopardizes biodiversity in rock--paper--scissors games. {\em Nature},
  \textbf{448}(7157), 1046--1049.

\bibitem[Reichenbach et~al., 2008]{reichenbach2008}
Reichenbach, T., Mobilia, M. {\&} Frey, E. (2008)  Self-organization of mobile
  populations in cyclic competition. {\em Journal of Theoretical Biology},
  \textbf{254}(2), 368--383.

\bibitem[Saad, 2011]{saad2011numerical}
Saad, Y. (2011) {\em Numerical methods for large eigenvalue problems: revised
  edition}.
Philadelphia, PA: SIAM.

\bibitem[Sandstede \& Scheel, 2000]{sandstede2000absolute}
Sandstede, B. {\&} Scheel, A. (2000)  Absolute versus convective instability of
  spiral waves. {\em Physical Review E}, \textbf{62}(6), 7708.

\bibitem[Sandstede \& Scheel, 2020]{SandScheel2020}
Sandstede, B. {\&} Scheel, A. (2020)  Spiral waves: linear and nonlinear
  theory. {\em arXiv preprint arXiv:2002.10352}.

\bibitem[Sinervo \& Lively, 1996]{Sinervo1996}
Sinervo, B. {\&} Lively, C.~M. (1996)  The rock--paper--scissors game and the
  evolution of alternative male strategies. {\em Nature}, \textbf{380}(6571),
  240--243.

\bibitem[Sinervo et~al., 2000]{Sinervo2000}
Sinervo, B., Miles, D.~B., Frankino, W.~A., Klukowski, M. {\&} DeNardo, D.~F.
  (2000)  Testosterone, endurance, and Darwinian fitness: natural and sexual
  selection on the physiological bases of alternative male behaviors in
  side-blotched lizards. {\em Hormones and Behavior}, \textbf{38}(4), 222--233.

\bibitem[Szczesny et~al., 2013]{Szczesny2013}
Szczesny, B., Mobilia, M. {\&} Rucklidge, A.~M. (2013)  When does cyclic
  dominance lead to stable spiral waves?. {\em EPL (Europhysics Letters)},
  \textbf{102}(2), 28012.

\bibitem[Szczesny et~al., 2014]{Szczesny2014}
Szczesny, B., Mobilia, M. {\&} Rucklidge, A.~M. (2014)  Characterization of
  spiraling patterns in spatial rock-paper-scissors games. {\em Physical Review
  E}, \textbf{90}(3), 032704.

\bibitem[Szolnoki et~al., 2014]{Szolnoki2014}
Szolnoki, A., Mobilia, M., Jiang, L.-L., Szczesny, B., Rucklidge, A.~M. {\&}
  Perc, M. (2014)  Cyclic dominance in evolutionary games: a review. {\em
  Journal of the Royal Society Interface}, \textbf{11}(100), 20140735.

\bibitem[Taylor \& Aarssen., 1990]{Taylor1990}
Taylor, D.~R. {\&} Aarssen., L.~W. (1990)  Complex competitive relationships
  among genotypes of three perennial grasses: implications for species
  coexistence. {\em The American Naturalist}, \textbf{163}(3), 305--327.

\bibitem[Woods \& Champneys, 1999]{woods1999heteroclinic}
Woods, P. {\&} Champneys, A. (1999)  Heteroclinic tangles and homoclinic
  snaking in the unfolding of a degenerate reversible Hamiltonian--Hopf
  bifurcation. {\em Physica D: Nonlinear Phenomena}, \textbf{129}(3-4),
  147--170.

\bibitem[Wulff \& Schebesch, 2006]{wulff2006}
Wulff, C. {\&} Schebesch, A. (2006)  Numerical continuation of symmetric
  periodic orbits. {\em SIAM Journal on Applied Dynamical Systems},
  \textbf{5}(3), 435--475.

\end{thebibliography}

\end{document}